\documentclass[reqno,12pt]{amsart}
\usepackage{amsfonts,amsmath,amsthm,amssymb,amscd}
\newcommand{\version}{Ver.~0.0}
\newcommand{\setversion}[1]{\renewcommand{\version}{Ver.~{#1}}}
\setversion{1.0 [Mon Jan 26 19:55:53 JST 2004]}
\setversion{1.1 [Tue Jan 27 18:39:18 JST 2004]}
\setversion{1.2 [Wed Jan 28 15:06:01 JST 2004]}
\setversion{1.3 [Fri Feb  6 16:58:52 MET 2004]}
\setversion{2.0 [Mon Feb  9 19:27:27 JST 2004]}
\setversion{2.1 [Fri Apr 23 19:11:40 JST 2004]}
\setversion{2.2 [Fri Jun 25 14:38:27 JST 2004]}

\title[Equivariant double fibrations]
{A note on affine quotients and equivariant double fibrations}
\thanks{\textit{Proceedings of ``Infinite Dimensional Harmonic Analysis {\upshape(}T\"{u}bingen, 2003/9/14 -- 9/21{\upshape)}''}}

\author{Kyo Nishiyama}
\address{
Department of Mathematics\\
Graduate School of Science\\
Kyoto University\\
Sakyo, Kyoto 606-8502, Japan}
\email{kyo@math.kyoto-u.ac.jp}


\theoremstyle{plain}
\newtheorem{theorem}{Theorem}

\newtheorem{corollary}[theorem]{Corollary}
\newtheorem{lemma}[theorem]{Lemma}

\newtheorem{assumption}[theorem]{Assumption}
\newtheorem{introtheorem}{Theorem}

\newtheorem{introassumption}[introtheorem]{Assumption}
\theoremstyle{definition}
\newtheorem{definition}[theorem]{Definition}
\newtheorem{example}[theorem]{\upshape Example}
\newtheorem{introproblem}[introtheorem]{Problem}

\theoremstyle{remark}

\numberwithin{equation}{section}
\numberwithin{theorem}{section}

\newcommand{\Z}{\mathbb{Z}}

\newcommand{\R}{\mathbb{R}}
\newcommand{\bbH}{\mathbb{H}}

\newcommand{\C}{\mathbb{C}}

\newcommand{\lie}[1]{\mathfrak{#1}}

\newcounter{thmenum}
\newenvironment{thmenumerate}{%
\begin{list}{{\upshape(\thethmenum)}}{%
\usecounter{thmenum}
\setlength{\labelsep}{.5em}
\setlength{\labelwidth}{0pt}
\setlength{\parsep}{0pt}
\setlength{\leftmargin}{3pt}
\setlength{\rightmargin}{0pt}
\setlength{\itemindent}{\leftmargin}
}}
{\end{list}}

\newcommand{\mycomment}[1]{} % Nothing to do

\newlength{\lengthcup}
\newcommand{\transpose}[1]{\,{}^t{#1}}

\newcommand{\Hom}{\mathop\mathrm{Hom}\nolimits{}}

\newcommand{\rank}{\mathop\mathrm{rank}\nolimits{}}

\newcommand{\directsum}{\mathop{\ \sum\nolimits^{\oplus}}}

\newcommand{\closure}[1]{\overline{#1}}

\newcommand{\trivial}{\mathbf{1}}

\newcommand{\partition}{\mathcal{P}}

\newcommand{\Spec}{\mathop\mathrm{Spec}\nolimits{}}
\newcommand{\Sym}{\mathrm{Sym}}
\newcommand{\Alt}{\mathrm{Alt}}

\newcommand{\orbit}{\mathbb{O}}  % 
\newcommand{\nullcone}{\lie{N}}

\newcommand{\harmonics}{\mathcal{H}}

\newcommand{\composit}{\odot}

\newcommand{\Gbig}{G}

\newcommand{\Gsmall}{G'}

\newcommand{\Obig}{\orbit}
\newcommand{\Osmall}{{\orbit}'}

\newcommand{\Obigtrivial}{\orbit{}^{\trivial}}

\newcommand{\GITquotient}{/\!/}

\newcommand{\AffineGrass}{\mathbb{G}^{\text{\upshape{aff}}}}
\newcommand{\Det}{\mathrm{Det}}
\newcommand{\calorbit}{\mathcal{O}}

\newcommand{\generic}[1]{{#1}^{\circ}}

\newcommand{\Djokovic}{{\hbox to.7pt{--\hss}D}okovi\'c}

\newcommand{\tbl}[2]{#2}

\begin{document}

\begin{abstract}
We consider two linear algebraic groups $ \Gbig $ and $ \Gsmall $ over the field of complex numbers, both of which are reductive.  
Take a finite dimensional rational representation $ W $ of $ \Gbig \times \Gsmall $.  
Let $ Y = W \GITquotient \Gbig := \Spec \C[W]^{\Gbig} $ and $ X = W \GITquotient \Gsmall := \Spec \C[W]^{\Gsmall} $ be the affine quotients.  
Since the action of $ \Gbig $ and $ \Gsmall $ commutes on $ W $, 
the quotient space $ X $ (respectively $ Y $) naturally inherits the action of $ \Gbig $ (respectively $ \Gsmall $).  

In this note, we study the interrelation between the orbit structures of $ X / \Gbig $ and $ Y / \Gsmall $.  
In a good situation, we can embed $ Y / \Gsmall $ into $ X / \Gbig $, and the embedding map $ \theta $ preserves important properties such as 
the closure relation and nilpotency.  
We give a sufficient condition for the existence of such embedding, and provide many examples 
arising from the natural representations of classical groups.  

As an application we consider the geometric problem of unimodular congruence classes of bilinear forms proposed by 
\Djokovic-Sekiguchi-Zhao.
\end{abstract}

\maketitle

\tableofcontents
%\newpage

\section*{Introduction}

Let us consider two linear algebraic groups $ \Gbig $ and $ \Gsmall $ over the field of complex numbers, both of which are reductive.  
Take a finite dimensional rational representation $ W $ of $ \Gbig \times \Gsmall $.  
The affine quotient of $ W $ by the action of $ \Gbig $ is denoted by $ Y = W \GITquotient \Gbig := \Spec \C[W]^{\Gbig} $, and 
similarly, $ X = W \GITquotient \Gsmall := \Spec \C[W]^{\Gsmall} $, where 
$ \C[W] $ is the ring of regular functions on $ W $, and the superscript of $ \Gbig $ denotes the subring of $ \Gbig $-invariants.  
Since the action of $ \Gbig $ and $ \Gsmall $ commutes on $ W $, 
the quotient space $ X $ (respectively $ Y $) naturally inherits the action of $ \Gbig $ (respectively $ \Gsmall $).  

In this note, we study the interrelation between the orbit structures of $ X / \Gbig $ and $ Y / \Gsmall $.  
In a good situation, we can embed $ Y / \Gsmall $ into $ X / \Gbig $, and the embedding map $ \theta $ preserves important properties such as 
the closure relation and nilpotency.  
We give a sufficient condition for the existence of such embedding.  
Let us briefly explain the condition.  

Let $ \varphi : W \to X $ and $ \psi : W \to Y $ be the quotient maps.  
Then $ \nullcone = \psi^{-1}( \psi( 0 ) ) $ is called the null cone (for the action of $ \Gbig $).  
We assume the pair $ ( \Gbig , \Gsmall ) $ and the representation $ W $ satisfy the following.

\begin{introassumption}
\label{assumption}
{\renewcommand{\thethmenum}{\alph{thmenum}}
\begin{thmenumerate}
\item
\label{assumption.1}
The quotient map $ \psi : W \rightarrow Y $ is flat.  This means the regular function ring $ \C[W] $ is flat over $ \C[W]^{\Gbig} $.  
\item
\label{assumption.2}
There exists an open dense $ \Gbig $-orbit $ \mathcal{O}_0 $ in $ \nullcone $.
\item
\label{assumption.3}
The null cone $ \nullcone $ is isomorphic to the scheme theoretic fiber $ W \times_Y \{ 0 \} $, 
i.e., the fiber product $ W \times_Y \{ 0 \} $ is reduced.  
\item
\label{assumption.4}
A generic fiber of 
the quotient map $ \varphi : W \rightarrow X $ is a single {\upshape(}hence closed{\upshape)} $ \Gsmall $-orbit.  
\item
\label{assumption.5}
Let $ \generic{W} $ be the union of closed $ \Gsmall $-orbits $ \calorbit' $ in $ W $ such that $ \varphi^{-1}( \varphi( \calorbit') ) = \calorbit' $.  
Then, for any $ y \in Y $, the fiber $ \psi^{-1}( y ) $ intersects $ \generic{W} $ non-trivially. 
\end{thmenumerate}
}
\end{introassumption}

Our main theorem, which is proved in \S \ref{sec:equivariant.fibration}, is the following.

\begin{introtheorem}
Let us assume Assumption {\upshape\ref{assumption}} holds.  
For any $ \Gsmall $-orbit $ \Osmall $ in $ Y $, there exists a $ \Gbig $-orbit $ \Obig $ in $ X $ such that $ \varphi( \psi^{-1}( \closure{\Osmall} ) ) = \closure{\Obig} $ holds.  
Thus we have a map $ \theta : Y / \Gsmall \to X / \Gbig $ which sends $ \Osmall $ to $ \Obig $.  
The lifting map $ \theta $ is injective and preserves the closure relation.  
If $ \Osmall \subset Y $ is a nilpotent orbit, then $ \Obig = \theta( \Osmall ) $ is also nilpotent.  
Moreover, we have $ \closure{\Obig} \simeq ( W \times_Y \closure{\Osmall} ) \GITquotient \Gsmall $.  
\end{introtheorem}

Let us summarize the brief history of the above theorem.  
The statement of the theorem is primitively noticed by early works of Roger Howe, 
and then clearly stated by Przebinda et al. for complex nilpotent orbits of reductive dual pairs 
(see e.g., \cite{Daszkiewicz.Przebinda.1996}, \cite{DKP.1997}).  
Recently, Daszkiewicz-Kra{\'s}kiewicz-Przebinda \cite{DKP.2002}, Ohta \cite{Ohta.preprint} 
and Nishiyama-Ochiai-Zhu \cite{NOZ.2003} extend the lifting map to the case of symmetric pairs arising from 
dual pairs of real reductive groups, but still only for the nilpotent orbits.  
By private communication (cf.~\cite{Ohta.Ohnuma}), 
T.~Ohta made me aware the fact that the orbits other than nilpotent ones are also in correspondence.  
In fact our proof of the lifting map in \cite{NOZ.2003} is applicable to all orbits without restriction, because of the geometric nature 
of the proof (see also \cite{Nishiyama.2000}).  
In this paper, we extend the correspondence to the general framework which is independent of the notion of dual pairs.

In \S \ref{sec:natural.representation}, we provide many examples which satisfy Assumption {\upshape\ref{assumption}}.  
The examples cover the cases in which $ W $ is the tensor product of the natural representations of classical groups, and 
also it contains several cases of contractive actions of the general linear groups.  
Although most of the cases are already obtained from the theory of dual pairs, 
two types of the lifting maps are newly found (see Theorems \ref{thm:tensor.product.type} and \ref{thm:contraction.by.GL}).

One of the new examples is strongly related to the $ SL(m, \C) $ action on the space of $ m \times m $-matrices $ M_m(\C) $; 
\begin{equation}
\label{eq:action.on.matrices}
A \mapsto g A \transpose{g} \qquad  ( g \in SL(m, \C) , \;  A \in M_m(\C) ) ,
\end{equation}
which is studied by 
D.~\v{Z}.~{\Djokovic}, 
J.~Sekiguchi and 
K.~Zhao \cite{DSZ.2003}; also the structure of the orbit space is being studied by H.~Ochiai recently.  

Our theory tells that the orbit space $ M_m(\C) / SL(m, \C) $ can be embedded into 
the orbit space of the affine cone of the Grassmann variety $ \AffineGrass_m( V \oplus V^{\ast} ) $ with $ GL(V) $-action.  

\begin{introtheorem}
There is an embedding map 
\begin{equation}
\theta : M_m(\C) / SL(m, \C) \to \AffineGrass_m( V \oplus V^{\ast} ) / GL(V) ,
\end{equation}
which preserves the closure relation, and carries nilpotent orbits to nilpotent ones.  
The image of the trivial orbit $ \Obigtrivial = \theta( \{ 0 \} ) $ is a spherical variety, and its closure in $ \AffineGrass_m( V \oplus V^{\ast} ) $ is normal.
\end{introtheorem}

In \S \ref{sec:SDZ.problem}, 
we have another embedding arising from the decomposition of the full matrix space into symmetric ones and skew-symmetric ones.  

Finally, we would like to propose some natural problems.  

\begin{introproblem}
\begin{thmenumerate}
\item
Find a pair $ ( \Gbig , \Gsmall ) $ and a representation $ W $ satisfying Assumption {\upshape\ref{assumption}}, 
for which one of the pair is an \textit{exceptional} group.  
\item
Consider irreducible representations $ V $ of $ \Gbig $ and $ U $ of $ \Gsmall $.  
Classify all the pairs $ ( V, U ) $ for which $ W = V \otimes U $ satisfies Assumption \ref{assumption}.
\item
Consider irreducible representations $ V $ of $ \Gbig = GL(n, \C) $ and $ U $ of $ \Gsmall $.  
Classify all the pairs $ ( V, U ) $ for which $ W = ( V \oplus V^{\ast} ) \otimes U $ satisfies Assumption \ref{assumption}.
\item
Give a complete description of the lifting map $ \theta $ in a combinatorial way.
\item
Find a representation theoretic interpretation of the lifting map $ \theta $.  
In the case of the liftings arising from dual pairs, it is provided by the \textit{theta correspondence} 
(or \textit{Howe correspondence}).  
See \cite{Daszkiewicz.Przebinda.1996} and \cite{NZ.2003}.
\item
Find the relation between the singularity of $ \closure{\Osmall} $ and that of $ \closure{\Obig} $, where $ \Obig = \theta( \Osmall) $ is the lifted orbit.
\end{thmenumerate}
\end{introproblem}

The author thanks Ralph Bremigan for useful discussion and for pointing out the
reference \cite{Schwarz.1978}.  

\section{Preliminaries}
\label{sec:preliminary}

In this section, we summarize definitions and well known facts on affine quotient maps.  

Let $ X $ be an affine variety on which a reductive algebraic group $ G $ acts rationally.  
We denote the affine coordinate ring (or ring of regular functions) on $ X $ by $ \C[X] $.  
Then $ G $ naturally acts on $ \C[X] $ via the formula
\begin{equation*}
g \cdot f(x) = f( g^{-1} \cdot x ) \qquad
( f(x) \in \C[X] , g \in G ) .
\end{equation*}
We denote the ring of $ G $-invariants in $ \C[X] $ by $ \C[X]^G $.  
The \textit{affine quotient} $ X \GITquotient G $ of $ X $ by the action of $ G $ is defined to be 
\begin{equation*}
X \GITquotient G = \Spec \C[X]^G .
\end{equation*}
The affine variety $ X \GITquotient G $ is often called the categorical quotient in the literature.  
The inclusion map $ \C[X]^G \hookrightarrow \C[X] $ induces a quotient morphism 
$ \varphi : X \to X \GITquotient G $, 
which has the following properties.

\begin{lemma}
Let $ \varphi : X \to X \GITquotient G $ be an affine quotient map.  
\begin{thmenumerate}
\item
For any $ y \in X \GITquotient G $, the fiber $ \varphi^{-1}( y ) $ is a $ G $-stable closed subvariety of $ X $, 
and it contains a unique closed $ G $-orbit.
\item
Let $ Z \subset X $ be a $ G $-stable closed subvariety.  
Then the restriction \linebreak[3]
$ \varphi \big|_Z : Z \longrightarrow \varphi(Z) \subset X \GITquotient G $ is an affine quotient map, and 
consequently $ \varphi(Z) \simeq Z \GITquotient G $.
\end{thmenumerate}
\end{lemma}

If we make a closed point $ y \in X \GITquotient G $ correspond to the unique closed $ G $-orbit in the fiber $ \varphi^{-1}( y ) \subset X $, 
we have a bijection between $ X \GITquotient G $ and the set of all closed $ G $-orbits in $ X $.  
In this sense, $ X \GITquotient G $ only classifies closed $ G $-orbits in $ X $.  

\medskip
In the following, we give three basic examples of affine quotient maps which will play important roles in the subsequent sections.

\subsection{Special linear group}
\label{subse:SL}

Let $ V = \C^n $ be a vector space on which $ G = SL(n, \C) $ acts naturally as the matrix multiplication.  
Take an another vector space $ U = \C^m $ and put $ W = V \otimes U $.  
$ G $ acts on $ W $ in the first component.  
If we identify $ W $ with the space of $ n \times m $-matrices over $ \C $, which we denote by $ M_{n,m} $, 
then the action is given by the matrix multiplication on the left.

Let us assume that $ n \leq m $ and identify $ W = M_{n,m} $.  
Then the ring of $ G $-invariants $ \C[W]^G $ is generated by all $ n \times n $-minors, 
which have the Pl\"{u}cker relations.  It is well known that there is no other relation among them 
(see, e.g., Theorem~3.1.6 in \cite{Manivel.2001}), 
and we can identify the quotient $ W \GITquotient G $ 
with the affine cone of the Grassmannian variety of $ n $-dimensional subspaces in $ U $.  
We denote it by $ \AffineGrass_n(U) $.  
The affine quotient map $ \varphi : W = M_{n, m} \to \AffineGrass_n(U) $ is interpreted as follows.  
By the Pl\"{u}cker embedding, we consider $ \AffineGrass_n(U) $ as the closed subvariety of $ \bigwedge^n U $.  
Under this identification, the quotient map $ \varphi $ sends $ A \in M_{n, m} $ to the exterior product of its rows.  
Thus, if $ \rank A < n $, then $ \varphi(A) = 0 \in \bigwedge^n U $.

If $ n > m $, then the only $ G $-invariants in $ \C[W] $ are scalars.  
So we have $ W \GITquotient G = \{ \ast \} $ (one point).

\subsection{General linear group}
\label{subsec:GL.quotient}

Let $ V = \C^n $ be a natural (or defining) representation of $ GL(n, \C) $.  
Take another vector spaces $ U^+ = \C^p $ and $ U^- = \C^q $, 
and put $ W = V \otimes U^+ \oplus V^{\ast} \otimes U^- $, where $ V^{\ast} $ denotes the contragredient representation of $ V $.  
$ G $ acts on $ W $ naturally in the first components.  
We identify $ W $ with $ \Hom( U^+{}^{\ast} , V ) \oplus \Hom( V, U^- ) $.

If $ n \geq p, q $, then it is easy to see that $ W \GITquotient G = U^+ \otimes U^- $.  
The quotient map $ \varphi $ is given by 
\begin{equation}
\label{eqn:quotient.map.GL}
\begin{split}
\varphi : W \simeq \Hom( U^+{}^{\ast} , V ) & \oplus \Hom( V, U^- ) \ni ( f , g ) 
\mapsto g \circ f \\
& \in \Hom( U^+{}^{\ast} , U^- ) \simeq U^+ \otimes U^- .  
\end{split}
\end{equation}

Let us consider the case where $ n < \max \{ p, q \} $.  
In this case, we have 
\begin{equation*}
W \GITquotient G = \Det_n( U^+ \otimes U^- ) , 
\end{equation*}
where $ \Det_n( U^+ \otimes U^- ) $ denotes 
the determinantal variety of rank $ n $, which is isomorphic to 
\begin{equation*}
\{ f \in \Hom( U^+{}^{\ast} , U^- ) \mid \rank f \leq n \}
\end{equation*}
under the identification $ \Hom( U^+{}^{\ast} , U^- ) \simeq U^+ \otimes U^- $.  
The quotient map is essentially the same as in \eqref{eqn:quotient.map.GL}.

\subsection{Quadratic space}
\label{subsec:quodratic.form}

Let $ V = \C^n $ be a vector space with a non-degenerate bilinear form, which we assume symmetric or skew-symmetric.  
We denote by $ G $ the group of isometries on $ V $ so that $ G $ is an orthogonal group $ O(n, \C) $ or a symplectic group $ Sp( n, \C) $ 
according as the form is symmetric or skew-symmetric.  
Note that $ n $ is necessarily even in the skew-symmetric case since the bilinear form is non-degenerate.  
Take an another vector space $ U = \C^m $ and put $ W = V \otimes U $.  

Let us first consider the symmetric case, hence $ G = O(n, \C) $.    

If $ n \geq m $, then the quotient $ W \GITquotient G $ is isomorphic to the symmetric tensor product $ \Sym(U) $ in $ U \otimes U $.  
Let us identify $ \Sym(U) $ with $ \{ h \in \Hom( U^{\ast} , U ) \mid \transpose{h} = h \} $, 
where $ \transpose{h} $ denotes the transpose of $ h $.  
For $ f \in \Hom( V^{\ast} , U ) \simeq W $, 
the image of the quotient map $ \varphi $ is given by 
\begin{equation*}
\varphi(f) : U^{\ast} \xrightarrow{\quad \transpose{f} \quad} V \simeq V^{\ast} \xrightarrow{\quad f \quad} U ,
\end{equation*}
where the isomorphism $ V \simeq V^{\ast} $ is induced by the symmetric bilinear form.  
It is easy to see that $ \varphi(f) $ belongs to $ \Sym(U) $.

If $ n < m $, the above image $ \varphi(f) $ belongs to 
\begin{equation*}
\Sym_n(U) := \Sym(U) \cap \Det_n( U \otimes U ) , 
\end{equation*}
and we have $ W \GITquotient G = \Sym_n(U) $.

The skew-symmetric case is similar.  
If $ n \geq m $, we have $ W \GITquotient G \simeq \Alt(U) $, where 
$ \Alt(U) $ denotes the set of skew-symmetric tensor products in $ U \otimes U $.  
Note that it is canonically isomorphic to $ \{ h \in \Hom( U^{\ast} , U ) \mid \transpose{h} = - h \} $.  
If $ n < m $, we have 
\begin{equation*}
W \GITquotient G = \Alt_n(U) := \Alt(U) \cap \Det_n( U \otimes U ) .
\end{equation*}

\section{Equivariant double fibration}
\label{sec:equivariant.fibration}

Let $ \Gbig $ and $ \Gsmall $ be connected linear algebraic groups over $ \C $ which are reductive.  
Suppose that there exists a finite dimensional complex vector space $ W $, on which $ \Gbig \times \Gsmall $ acts linearly.  
We put 
\begin{align*}
X &= W \GITquotient \Gsmall , \quad \text{with quotient map } \varphi : W \rightarrow X , \\
Y &= W \GITquotient \Gbig   , \quad \text{with quotient map } \psi : W \rightarrow Y .
\end{align*}
Then $ \Gbig $ naturally acts on $ X $, and similarly, $ Y $ inherits an action of $ \Gsmall $.  
By abuse of notation, we denote the image $ \varphi(0) $ (respectively $ \psi(0) $) of $ 0 \in W $ 
simply by $ 0 \in X $ (respectively $ 0 \in Y $).  

\begin{definition}
A $ \Gbig $-orbit $ \Obig \subset X $ is said to be \textit{nilpotent} 
if $ \closure{\Obig} $ contains $ 0 $.  
The same definition applies to a $ \Gsmall $-orbit $ \Osmall \subset Y $.
\end{definition}

Let $ Z = W \GITquotient (\Gbig \times \Gsmall) $ be the affine quotient of $ W $ by $ \Gbig \times \Gsmall $, 
which is naturally identified with $ X \GITquotient \Gbig $ and $ Y \GITquotient \Gsmall $ respectively.  
We denote the induced quotient maps by 
$ \psi_0 : X \rightarrow Z \simeq X \GITquotient \Gbig $ and 
$ \varphi_0 : Y \rightarrow Z \simeq Y \GITquotient \Gsmall $.  

\begin{lemma}
\label{lemma:preservation.of.nilpotent.orbit}
For a nilpotent $ \Gsmall $-orbit $ \Osmall \subset Y $, the subset 
$ \varphi( \psi^{-1} ( \closure{\Osmall} ) ) \subset X $ is a union of nilpotent $ \Gbig $-orbits.
\end{lemma}

\begin{proof}
Let $ \Obig \subset X $ be a $ \Gbig $-orbit. 
Then $ \Obig $ is nilpotent if and only if $ \psi_0( \closure{\Obig} ) = \{ 0 \} $, where $ 0 \in Z $ is the image of $ 0 \in W $.  
Thus, it is enough to show that the image of $ \varphi( \psi^{-1} ( \closure{\Osmall} ) ) $ under the map $ \psi_0 $ is $ \{ 0 \} $.  
Since $ \psi_0 \circ \varphi = \varphi_0 \circ \psi $, we have 
\begin{equation*}
\psi_0 \circ \varphi( \psi^{-1} ( \closure{\Osmall} ) ) = \varphi_0 \circ \psi ( \psi^{-1} ( \closure{\Osmall} ) ) 
= \varphi_0( \closure{\Osmall} ) = \{ 0 \} .
\end{equation*}
\end{proof}

Let $ \nullcone = \psi^{-1}( 0 ) \subset W $ be a \textit{null cone} (or null fiber).  
Throughout this article, we assume the following.

\begin{assumption}
\label{basic.assumption.Y}
{\renewcommand{\thethmenum}{\alph{thmenum}}
\begin{thmenumerate}
\item
\label{basic.assumption.1}
The quotient map $ \psi : W \rightarrow Y $ is flat.  This means the regular function ring $ \C[W] $ is flat over $ \C[W]^{\Gbig} $.  
\item
\label{basic.assumption.2}
There exists an open dense $ \Gbig $-orbit $ \calorbit_0 $ in $ \nullcone $.
\item
\label{basic.assumption.3}
The null cone $ \nullcone $ is isomorphic to the scheme theoretic fiber $ W \times_Y \{ 0 \} $, 
i.e., the fiber product $ W \times_Y \{ 0 \} $ is reduced.  
\end{thmenumerate}
}
\end{assumption}

Few remarks are in order.  
If the action of $ G $ on $ W $ is cofree (i.e., $ \C[W] $ is a graded free module over $ \C[W]^G $), then $ \psi $ is flat.  
The cofree actions are classified by J.~Schwarz \cite{Schwarz.1978}.  
The assumptions \eqref{basic.assumption.2} and \eqref{basic.assumption.3} imply that 
the null cone $ \nullcone \simeq W \times_Y \{ 0 \} $ is reduced and irreducible.
The irreducibility follows from the assumption \eqref{basic.assumption.2}.  
Moreover, if $ \Gbig $ is semisimple, the assumption \eqref{basic.assumption.2} implies that $ W \times_Y \{ 0 \} $ is reduced 
(see Korollar~2 in \cite{Knop.1986}), 
hence \eqref{basic.assumption.3} holds automatically.  

\medskip
Under these assumptions, we have

\begin{theorem}
\label{thm:lifting.of.orbit}
Take a $ \Gsmall $-orbit $ \Osmall $ in $ Y $.  
\begin{thmenumerate}
\item
The scheme theoretic inverse image $ \psi^{-1}( \closure{\Osmall} ) = W \times_Y \closure{ \Osmall } $ is reduced and irreducible.
\item
The inverse image $ \psi^{-1}( \closure{ \Osmall } ) $ contains an open dense $ \Gbig \times \Gsmall $-orbit $ \mathfrak{O} $, hence 
there is a $ \Gbig $-orbit $ \Obig $ in $ X $ such that $ \varphi( \psi^{-1}( \closure{ \Osmall } ) ) = \closure{ \Obig } $.
We say the $ \Gbig $-orbit $ \Obig $ is \textit{lifted} from $ \Osmall $, and denote it by $ \Obig = \theta( \Osmall ) $.
\item
The lifting map $ \theta $ preserves the closure relation.  
If $ \Osmall $ is a nilpotent $ \Gsmall $-orbit, then $ \Obig = \theta( \Osmall ) $ is also nilpotent.
\end{thmenumerate}
\end{theorem}

\begin{proof}
This theorem is a generalization of Theorems 2.5 and 2.10 in \cite{NOZ.2003}.  
Note that the results in \cite{NOZ.2003} are stated for nilpotent orbits, but actually they are valid for all kind of orbits.  
Hence, the proof is almost the same as in \cite{NOZ.2003}, 
but for convenience of the reader, we indicate the outline of the proof.  

First, we prove that the scheme theoretic fiber $ W \times_Y \{ y \} $ is reduced for any $ y \in Y $.  
Then this will imply that $ W \times_Y Z $ is reduced for an arbitrary closed subvariety $ Z \subset Y $.  

In the terminology of commutative algebra, the claim that $ W \times_Y \{ y \} $ is reduced is equivalent to 
that $ \C[W] \otimes_{\C[Y]} \C_y $ does not contain any non-zero nilpotent element, where $ \C_y $ denotes the function ring on 
the one point set $ \{ y \} $.  
Note that we assume that $ \C[W] \otimes_{\C[Y]} \C_0 \simeq \C[\nullcone] $ is an integral domain.  
Since $ \C[W] \otimes_{\C[Y]} \C_y $ is a deformation of the homogeneous integral domain $ \C[\nullcone] $, it is also an integral domain.  
Thus we have proven $ W \times_Y \{ y \} $ is reduced and irreducible.

Next, we shall prove that the fiber $ \psi^{-1}( y ) $ contains a dense open $ \Gbig $-orbit.  
Put $ M = \psi^{-1}(y) $ and denote by $ \widehat{M} $ the asymptotic cone of $ M $ 
(see \S~5.2 of \cite{Popov.Vinberg.1994} for the definition of asymptotic cones).  
Then, by the flatness of $ \psi $, the asymptotic cone $ \widehat{M} $ coincides with the null cone $ \nullcone $.  
Let $ \calorbit_y $ be a generic $ \Gbig $-orbit in $ M $.  
Consider the cone $ \C M $ generated by $ M $ in 
$ W $, then it is clear that the dimension of a generic orbit in $ \C M $ is equal to $ \dim \calorbit_y $, which  
in turn coincides with the dimension of the generic orbit in $ \closure{\C M} \subset W $.  
Since $ \nullcone = \widehat{M} \subset \closure{\C M} $, 
the dimension of a generic orbit in $ \nullcone $ cannot exceed that of $ \calorbit_y $.  
Note that $ \nullcone $ has an open dense orbit by Assumption \ref{basic.assumption.Y}~\eqref{basic.assumption.2}.
This means that $ \dim \calorbit_y \geq \dim \nullcone $.  
On the other hand, we have the equality $ \dim M = \dim \widehat{M} = \dim \nullcone $ of dimensions, 
hence $ \dim \calorbit_y \geq \dim M $.  
Since $ \calorbit_y \subset M $, we conclude that $ \dim \calorbit_y = \dim M $, and 
that $ \calorbit_y $ is an open dense orbit in $ M $, by the irreducibility of $ M $ just proved above.  

Since $ \psi $ is $ \Gsmall $-equivariant, we get 
$ \psi^{-1}( \Osmall ) = \Gsmall \cdot \psi^{-1}( y ) $ for any $ y \in \Osmall $.  
Note that $ \psi^{-1}( y ) $ contains an open dense $ \Gbig $-orbit $ \calorbit_y $.  
Choose an arbitrary point $ w \in \calorbit_y $, and we see the $ \Gbig \times \Gsmall $-orbit 
$ \mathfrak{O} = \Gsmall \Gbig w $ is open dense in $ \closure{ \psi^{-1}( \Osmall ) } $.  

Since we assume that $ \psi $ is flat, it is an open map by Ex.~(III.9.1) in \cite{Hartshorne.1977}.  
Thus the equality $ \closure{ \psi^{-1}( \Osmall ) } = \psi^{-1}( \closure{ \Osmall } ) $ holds.  
Now we conclude that $ \Obig = \Gbig \varphi(w) \subset \varphi( \psi^{-1}( \closure{ \Osmall } ) ) $ is the open dense orbit which we want.

The claim that the lifting map preserves the closure relation is obvious from the definition of $ \theta $.  
Lemma \ref{lemma:preservation.of.nilpotent.orbit} tells us that $ \theta $ preserves nilpotency.  
\end{proof}

\begin{corollary}
Let $ \Osmall $ be a  $ \Gsmall $-orbit in $ Y $, and $ \Obig = \theta( \Osmall ) $ its lift.  
Then we have
\begin{equation*}
\C[ \closure{ \Obig } ] \simeq \Bigl( \C[ W ] \otimes_{\C[Y]} \C[ \closure{\Osmall} ] \Bigr){}^{\Gsmall} .
\end{equation*}
\end{corollary}

If $ \C[W] $ is free over $ \C[W]^{\Gsmall} $, we can write 
$ \C[W] = \harmonics \otimes \C[W]^{\Gsmall} $, where $ \harmonics $ is the space of $ \Gsmall $-harmonic polynomials in $ \C[W] $.  
Then, the above corollary tells us that 
\begin{equation*}
\C[ \closure{ \Obig } ] \simeq \bigl( \harmonics \otimes \C[ \closure{\Osmall} ] \bigr){}^{\Gsmall} .
\end{equation*}
Note that, as a $ \Gsmall $-module, $ \harmonics $ is isomorphic to the regular function ring $ \C[ \nullcone ] $ of $ \nullcone $.  

Let us denote by $ X / \Gbig $ the set of all $ \Gbig $-orbits.  
Note that $ X / \Gbig $ may not be an algebraic variety, but only a topological space.  

In general, the lifting map $ \theta : Y / \Gsmall \to X / \Gbig $ is not necessarily injective.  
Let us give a sufficient condition for the injectivity of $ \theta $.  We denote
\begin{equation}
\label{eqn:definition.of.generic.orbits}
\begin{split}
\generic{X} &= \{ x \in X \mid \varphi^{-1}(x) \text{ consists of a single $ \Gsmall $-orbit} \} , \text{ and } \\
\generic{W} &= \varphi^{-1}( \generic{X} ) = \coprod_{x \in \generic{X}} \varphi^{-1}( x ) .
\end{split}
\end{equation}

\begin{theorem}
\label{thm:injectivity}
If, for any $ y \in Y $, the fiber $ \psi^{-1}( y ) $ intersects $ \generic{W} $ non-trivially,  
then the lifting map $ \theta : Y / \Gsmall \to X / \Gbig $ is injective.
\end{theorem}

\begin{proof}
Let $ \Osmall_1 \neq \Osmall_2 $ be two different $ \Gsmall $-orbits in $ Y $, and denote $ \Obig_i = \theta( \Osmall_i ) \subset X \; ( i = 1, 2 ) $.  
Without loss of generality, we can assume that $ \Osmall_1 \cap \closure{\Osmall_2} = \emptyset $.  
Then $ \psi^{-1}( \Osmall_1 ) \cap \psi^{-1}( \closure{\Osmall_2} ) = \emptyset $ and, by the assumption, 
$ \psi^{-1}( \Osmall_1 ) $ contains a closed $ \Gsmall $-orbit which is of the form $ \varphi^{-1}( x ) $ for some $ x \in X $.  
This means that $ x \not\in \varphi( \psi^{-1}( \closure{\Osmall_2} ) ) = \closure{\Obig_2} $, 
while $ x \in \varphi( \psi^{-1}( \closure{\Osmall_1} ) ) = \closure{\Obig_1} $.  
Thus $ \Obig_1 \neq \Obig_2 $ which proves the theorem.  
\end{proof}

\section{Double fibration related to the natural representations}
\label{sec:natural.representation}

Here we give several examples which satisfy Assumption \ref{basic.assumption.Y}.  
These examples arise from the natural representations of various classical groups.  

To exclude trivial cases, we further assume the following.  

\begin{assumption}
\label{basic.assumption.X}
{\renewcommand{\thethmenum}{\alph{thmenum}}
\begin{thmenumerate}
\setcounter{thmenum}{3}
\item
\label{basic.assumption.X:item1}
A generic fiber of 
the quotient map $ \varphi : W \rightarrow X $ is a single {\upshape(}hence closed{\upshape)} $ \Gsmall $-orbit.  
\item
\label{basic.assumption.X:item2}
Put $ \generic{X} $ and $ \generic{W} $ as in \eqref{eqn:definition.of.generic.orbits}.  
For any $ y \in Y $, the fiber $ \psi^{-1}( y ) $ intersects $ \generic{W} $ non-trivially. 
\end{thmenumerate}
}
\end{assumption}

Assumptions \ref{basic.assumption.Y} and \ref{basic.assumption.X} assure that the lifting map 
$ \theta : Y / \Gsmall \to X / \Gbig $ is injective, and preserves the closure ordering.  
Also $ \theta $ lifts nilpotent orbits to nilpotent orbits.

\subsection{Tensor products}  

We first investigate examples satisfying Assumption \ref{basic.assumption.Y}.  
Let $ V $ be a finite dimensional representation of $ \Gbig $, and $ U $ a finite dimensional vector space.  
Put $ W = V \otimes U $ on which $ \Gbig $ acts naturally.  

\begin{lemma}
The quotient map $ \psi : W = V \otimes U \rightarrow Y := W \GITquotient \Gbig $ satisfies Assumption {\upshape\ref{basic.assumption.Y}} if 
the representation $ ( \Gbig , V ) $ and a vector space $ U $ are in Table~{\upshape\ref{table:1}}.
Here we denote by $ \Sym(U) $ {\upshape (}respectively $ \Alt(U) ${\upshape )} the set of symmetric 
{\upshape (}respectively alternating{\upshape )} tensor products in $ U \otimes U $.
In these cases, the action of $ G $ on $ W $ is cofree.
\begin{table}[hbtp]
\caption{}
\tbl{}{
%\begin{center}
\setlength{\arraycolsep}{7pt}
\renewcommand{\arraystretch}{1.1}
\begin{math}\displaystyle
\begin{array}{c|c|c|c}
\Gbig & V  & U  & Y \\ \hline
\\[-8pt]
O(n, \C)  & \C^n \; \text{\upshape (natural)} & 2 \dim U < n & \Sym(U) \\
Sp(2 n, \C) & \C^{2 n} \; \text{\upshape (natural)} & \dim U \leq n & \Alt(U) 
\end{array}
\end{math}
%\end{center}
\label{table:1}
}
\end{table}
\end{lemma}

Next, let us consider examples satisfying Assumption \ref{basic.assumption.X}~\eqref{basic.assumption.X:item1}, 
i.e., we need to check that a generic fiber of the quotient map is a single orbit.  
Let $ U $ be a finite dimensional representation of $ \Gsmall $.  
Take an arbitrary finite dimensional vector space $ V $ and put $ W = V \otimes U $ as above.   

\begin{lemma}
The quotient map $ \varphi : W = V \otimes U \rightarrow X = W \GITquotient \Gsmall $ satisfies 
Assumption {\upshape\ref{basic.assumption.X}~\eqref{basic.assumption.X:item1}} if 
the representation $ ( \Gsmall , U ) $ and a vector space $ V $ are in
Table~{\upshape\ref{table:2}}.
For the notation of $ \Sym_m(V), \Alt_m(V) $ and $ \AffineGrass_m(V) $, see {\upshape\S \ref{subse:SL}} and 
{\upshape\S \ref{subsec:quodratic.form}}.
\begin{table}[hbtp]
\caption{}
\tbl{}{
\setlength{\arraycolsep}{7pt}
\renewcommand{\arraystretch}{1.1}
%\begin{center}
\begin{math}\displaystyle
\begin{array}{c|c|c|c}
\Gsmall & U  & V  & X \\ \hline
\\[-10pt]
O(m, \C)  & \C^m \; \text{\upshape (natural)} & \dim V \geq 1 & \Sym_m(V) \\
Sp(2 m, \C) & \C^{2 m} \; \text{\upshape (natural)} & \dim V \geq 1 & \Alt_m(V) \\
SL(m, \C) & \C^m  \; \text{\upshape (natural)} & \dim V \geq m & \AffineGrass_m(V)
\end{array}
\end{math}
%\end{center}
\label{table:2}
}
\end{table}
\end{lemma}

Combining these lemmas, we have the following 

\begin{theorem}
\label{thm:tensor.product.type}
Let $ W = V \otimes U $ be a representation of $ \Gbig \times \Gsmall
$ in Table~{\upshape\ref{table:3}}.
Here, for example, 
$ O(n, \C) \otimes O(m, \C) $ means the tensor product of the natural representations of $ O(n, \C) $ and $ O(m, \C) $.  
\begin{table}[hbtp]
\caption{}
\tbl{}{
\setlength{\arraycolsep}{3.5pt}
\renewcommand{\arraystretch}{1.1}
%\begin{center}
\begin{math}\displaystyle
\begin{array}{c|c|cc|c}
\qquad W & \text{\upshape condition } & X & Y & \text{\upshape dual pair} \\ \hline
&&&&\\[-10pt]
O(n, \C) \otimes O(m, \C) & 2 m < n & \Sym_m(\C^n) & \Sym(\C^m) & ( GL(n, \R) , GL(m, \R) ) \\
O(n, \C) \otimes Sp(2 m, \C) & 4 m < n & \Alt_{2 m}(\C^n) & \Sym(\C^{2 m}) & ( O(n, \C), Sp(2 m, \C) ) \\
O(n, \C) \otimes SL(m, \C) & 2 m < n & \AffineGrass_m(\C^n) & \Sym(\C^m) & \text{\upshape none} \\
Sp(2 n, \C) \otimes O(m, \C) & m \leq n & \Sym_m(\C^{2 n}) & \Alt(\C^m) & ( Sp(2 n, \C ) , O(m, \C) ) \\
Sp(2 n, \C) \otimes Sp(2 m, \C) & 2 m \leq n & \Alt_{2 m}(\C^{2 n}) & \Alt(\C^{2 m}) & ( GL(n, \bbH ) , GL( m, \bbH ) ) \\
Sp(2 n, \C) \otimes SL(m, \C) & m \leq n & \AffineGrass_m(\C^{2 n}) & \Alt(\C^m) & \text{\upshape none} \\
\end{array}
\end{math}
%\end{center}
\label{table:3}
}
\end{table}
\par
Then the double fibration by the affine quotient maps 
\begin{equation*}
X = W \GITquotient \Gsmall \xleftarrow{ \quad \varphi \quad } W \xrightarrow{ \quad \psi \quad } W \GITquotient \Gbig = Y
\end{equation*}
satisfies Assumptions {\upshape\ref{basic.assumption.Y}} and {\upshape\ref{basic.assumption.X}}.  
In particular, a  $ \Gsmall $-orbit $ \Osmall \subset Y $ is lifted to a  $ \Gbig $-orbit $ \Obig \subset X $, 
and the lifting map $ \theta $ is injective.
\end{theorem}

\subsection{Contraction by the action of a general linear group}

In this subsection, we consider the action of general linear groups.  

\begin{lemma}
Let $ V = \C^n $ be the natural {\upshape(}or defining{\upshape)} representation of $ \Gbig = GL(n, \C) $, and $ U^{\pm} $ finite dimensional vector spaces.  
Put $ W = V \otimes U^+ \oplus V^{\ast} \otimes U^- $ on which $ \Gbig $ acts naturally.  
Then the quotient map $ \psi : W \rightarrow Y = W \GITquotient \Gbig $ satisfies Assumption {\upshape\ref{basic.assumption.Y}} if and only if 
$ \dim U^+ + \dim U^- \leq \dim V $.  
In this case, the action of $ G $ on $ W $ is cofree   
and $ Y $ is naturally identified with $ U^+ \otimes U^- $.
\end{lemma}

\begin{lemma}
Let $ U = \C^m $ be the natural representation of $ \Gsmall = GL(m, \C) $, and $ V^{\pm} $ finite dimensional vector spaces.  
Put $ W = V^+ \otimes U \oplus V^- \otimes U^{\ast} $ on which $ \Gsmall $ acts naturally.  
Then the quotient map $ \varphi : W \rightarrow X = W \GITquotient \Gsmall $ satisfies 
Assumption {\upshape\ref{basic.assumption.X}~\eqref{basic.assumption.X:item1}} if and only if 
$ \dim V^+ = \dim V^- $ or 
$ \dim V^{\pm} \geq \dim U $.  
The quotient space $ X $ is naturally identified with the determinantal variety $ \Det_m( V^+ \otimes V^-) $ of rank 
$ m = \dim U $ {\upshape(}see \S {\upshape\ref{subsec:GL.quotient})}.  
\end{lemma}

\begin{theorem}
\label{thm:contraction.by.GL}
If $ W $ is one of the representations of $ \Gbig \times \Gsmall $
which are listed in \eqref{thm:contraction.by.GL:item1}--\eqref{thm:contraction.by.GL:item6} below, then the quotient maps
\begin{equation*}
X = W \GITquotient \Gsmall \xleftarrow{ \quad \varphi \quad } W \xrightarrow{ \quad \psi \quad } W \GITquotient \Gbig = Y
\end{equation*}
satisfy Assumptions {\upshape\ref{basic.assumption.Y}} and {\upshape\ref{basic.assumption.X}}.  
In particular, a  $ \Gsmall $-orbit $ \Osmall \subset Y $ is lifted to a  $ \Gbig $-orbit $ \Obig \subset X $, 
and the lifting map $ \theta $ is injective.
\begin{thmenumerate}
\item
\label{thm:contraction.by.GL:item1}
Let $ \Gbig = GL(n, \C) $ and $ V = \C^n $ the natural representation of $ \Gbig $.  
We put $ W = ( V \oplus V^{\ast} ) \otimes U $ for the natural 
representation $ U $ of $ \Gsmall $ in Table~{\upshape\ref{table:4}}.  
The quotient space $ X $ is given in Table~{\upshape\ref{table:4}} and $ Y = U \otimes U $.  
\begin{table}[hbtp]
\caption{}
\tbl{}{
\setlength{\arraycolsep}{7pt}
\renewcommand{\arraystretch}{1.1}
%\begin{center}
\begin{math}\displaystyle
\begin{array}{c|c|c}
\quad \Gsmall & \text{\upshape condition } &  \quad X \\ \hline
&&\\[-10pt]
O(m, \C) & 2 m \leq n & \Sym_m( V \oplus V^{\ast} ) \\
Sp(2 m, \C) & 4 m \leq n & \Alt_m( V \oplus V^{\ast} ) \\
SL(m, \C) & 2 m \leq n & \AffineGrass_m( V \oplus V^{\ast} )  
\end{array}
\end{math}
%\end{center}
\label{table:4}
}
\end{table}
\item
\label{thm:contraction.by.GL:item2}
Let $ \Gbig = GL(n, \C) $ and $ V = \C^n $ the natural representation of $ \Gbig $.  
We put $ \Gsmall = \Gsmall_+ \times \Gsmall_- $,   
and $ W = V \otimes U^+ \oplus V^{\ast} \otimes U^- $ for the natural 
representations $ U^{\pm} $ of $ \Gsmall_{\pm} $ in Table~{\upshape\ref{table:5}}.  
The quotient space $ X $ is given in Table~{\upshape\ref{table:5}} and $ Y = U^+ \otimes U^- $.  
\begin{table}[hbtp]
\caption{}
\tbl{}{
\setlength{\arraycolsep}{4pt}
\renewcommand{\arraystretch}{1.1}
%\begin{center}
\begin{math}\displaystyle
\begin{array}{cc|c|c|c}
\Gsmall_+ & \Gsmall_- & \text{\upshape condition } & X & \text{\upshape dual pair} \\ \hline
&&&&\\[-10pt]
O(p, \C) & O(q, \C) & p + q \leq n & \Sym_p(V) \oplus \Sym_q(V^{\ast}) & ( O( p, q ), Sp(2 n, \R) ) \\
Sp(2p, \C) & Sp(2 q, \C) & 2 p + 2 q \leq n & \Alt_{2 p}(V) \oplus \Alt_{2 q}(V^{\ast}) & ( Sp( 2 p  , 2 q ), O^{\ast}(2 n) ) \\
O(p, \C) & Sp( 2 q , \C ) & p + 2 q \leq n & \Sym_p(V) \oplus \Alt_{2 q}(V^{\ast}) & \text{\upshape none} \\
SL(p, \C) & SL(q, \C) & p + q \leq n & \AffineGrass_p(V) \times \AffineGrass_q(V^{\ast}) & \text{\upshape none} 
\end{array}
\end{math}
%\end{center}
\label{table:5}
}
\end{table}
\item
\label{thm:contraction.by.GL:item3}
Let $ \Gsmall = GL(m, \C) $ and $ U = \C^m $ the natural representation of $ \Gsmall $.
We put $ W = V \otimes ( U \oplus U^{\ast} ) $ for the natural 
representation $ V $ of $ \Gbig $ in Table~{\upshape\ref{table:6}}.  
The quotient space $ Y $ is given in Table~{\upshape\ref{table:6}} and $ X = \Det_m( V \otimes V ) $.  
\begin{table}[hbtp]
\caption{}
\tbl{}{
\setlength{\arraycolsep}{7pt}
\renewcommand{\arraystretch}{1.1}
%\begin{center}
\begin{math}\displaystyle
\begin{array}{c|c|c}
\Gbig & \text{\upshape condition } &  Y \\ \hline
\\[-10pt]
O(n, \C) & 4 m < n & \Sym( U \oplus U^{\ast} ) \\
Sp(2 n, \C) & 2 m \leq n & \Alt( U \oplus U^{\ast} ) 
\end{array}
\end{math}
%\end{center}
\label{table:6}
}
\end{table}
\item
\label{thm:contraction.by.GL:item4}
Let $ \Gsmall = GL(m, \C) $ and $ U = \C^m $ the natural representation of $ \Gsmall $.  
We put $ \Gbig = \Gbig_+ \times \Gbig_- $,   
and $ W = V^+ \otimes U \oplus V^- \otimes U^{\ast} $ for the natural 
representations $ V^{\pm} $ of $ \Gbig_{\pm} $ in Table~{\upshape\ref{table:7}}.  
The quotient space $ Y $ is given in Table~{\upshape\ref{table:7}} 
and $ X = \Det_m( V^+ \otimes V^- ) $.  
\begin{table}[hbtp]
\caption{}
\tbl{}{
\setlength{\arraycolsep}{4pt}
\renewcommand{\arraystretch}{1.1}
%\begin{center}
\begin{math}\displaystyle
\begin{array}{cc|c|c|c}
\Gbig_+ & \Gbig_- & \text{\upshape condition } & Y & \text{\upshape dual pair} \\ \hline
&&&&\\[-10pt]
O(r, \C) & O(s, \C) & 2 m < r, s & \Sym(U) \oplus \Sym(U^{\ast}) & ( Sp(2 m, \R) , O( r, s ) ) \\
Sp(2r, \C) & Sp(2 s, \C) & m \leq r, s &  \Alt(U) \oplus \Alt(U^{\ast}) & ( O^{\ast}(2 m) , Sp( 2r, 2s ) ) \\
O(r, \C) & Sp( 2 s , \C ) & 2 m \leq r - 1, s & \Sym(U) \oplus \Alt(U^{\ast}) & \text{\upshape none} 
\end{array}
\end{math}
%\end{center}
\label{table:7}
}
\end{table}
\item
\label{thm:contraction.by.GL:item5}
Let $ V = \C^n $ {\upshape(}respectively $ U = \C^m ${\upshape)}
be the natural representation of $ \Gbig = GL(n, \C) $
{\upshape(}respectively $ \Gsmall = GL(m, \C) ${\upshape)}  
and assume $ 2 m \leq n $.  
We put $ W = V \otimes U^{\ast} \oplus V^{\ast} \otimes U $.  
The quotient spaces are given by 
\begin{equation*}
X = \Det_m( V^{\ast} \otimes V) \quad \text{and} \quad Y = U^{\ast} \otimes U .  
\end{equation*}
The corresponding dual pair is $ ( GL(m, \C) , GL(n, \C) ) $.  
\item
\label{thm:contraction.by.GL:item6}
Let 
$ W =  
\C^p \otimes \C^r \oplus (\C^p)^{\ast} \otimes ( \C^s )^{\ast} \oplus 
( \C^q )^{\ast} \otimes ( \C^r )^{\ast} \oplus \C^q \otimes \C^s  $ 
be a representation of $ \Gbig \times \Gsmall $ with $ p + q \leq r, s $, 
where $ \Gbig = GL(r, \C) \times GL(s, \C) $ and $ \Gsmall = GL( p, \C ) \times GL( q , \C ) $.  
The quotient spaces are given by 
\begin{align*}
X &= \Det_p( \C^r \otimes ( \C^s )^{\ast} ) \times \Det_q( (\C^r)^{\ast} \otimes \C^s ) \quad \text{ and } \\
Y &= \C^p \otimes ( \C^q )^{\ast} \oplus (\C^p)^{\ast} \otimes \C^q .  
\end{align*}
The corresponding dual pair is $ ( U(p, q), U(r, s) ) $.  
\end{thmenumerate}
\end{theorem}

%\subsection{Remark on the lifting of orbits for symmetric pairs}

\section{Application to \Djokovic-Sekiguchi-Zhao problem}
\label{sec:SDZ.problem}

Recently, from the view point of the invariant theory, 
D.~\v{Z}.~{\Djokovic}, J.~Sekiguchi and K.~Zhao \cite{DSZ.2003} are studying the 
$ SL(m, \C) $-action on $ M_m(\C) = \C^m \otimes \C^m $ (see \eqref{eq:action.on.matrices}).  
The orbit structure of the action is also being studied by H.~Ochiai.
Here we resolve this action into two different ways, which 
fit our theory of the double fibration.

\subsection{Resolution via the contraction by the action of $ GL(n, \C) $.}
\label{subsec:resolution.by.contraction}

Let  $ \Gbig = GL(n, \C) $ and $ V = \C^n $ the natural representation of $ \Gbig $.  
We put $ W = ( V \oplus V^{\ast} ) \otimes U $ for the natural representation $ U = \C^m $ of $ \Gsmall = SL(m, \C) $.  
In the following, we assume that $ 2 m \leq n $. 

Theorem \ref{thm:contraction.by.GL} \eqref{thm:contraction.by.GL:item1} tells us that  the double fibration by the affine quotient maps 
\begin{equation*}
X = W \GITquotient \Gsmall \xleftarrow{ \quad \varphi \quad } W \xrightarrow{ \quad \psi \quad } W \GITquotient \Gbig = Y
\end{equation*}
satisfies Assumptions \ref{basic.assumption.Y} and \ref{basic.assumption.X}.  
Here we only check Assumption \ref{basic.assumption.X}~\eqref{basic.assumption.X:item2}.

\begin{lemma}
For any $ y \in Y $, the fiber $ \psi^{-1}( y ) $ intersects a closed $ \Gsmall $-orbit, 
which is precisely a fiber $ \varphi^{-1}( x ) $ for some $ x \in X $.
\end{lemma}

\begin{proof}
We identify $ W = V \otimes U \oplus V^{\ast} \otimes U $ with the space of $ 2 n \times m $-matrices $ M_{2n, m} $.  
Then a non-zero $ w \in M_{2n, m} $ generates a closed $ \Gsmall $-orbit if and only if $ \rank w = m $.  
In fact, for any non-zero $ x \in \AffineGrass_m( U ) $, 
we can prove that the fiber $ \varphi^{-1}( x ) $ consists of a single $ \Gsmall $-orbit, 
hence it is closed.  
Thus all non-closed orbits are contained in the null fiber $ \varphi^{-1}( 0 ) $ which is characterized by $ \rank w < m $.  

If we write $ w = ( A, B ) \in M_{n, m} \times M_{n, m} $, the quotient map $ \psi $ is given by 
\begin{equation*}
\psi(w) = \transpose{A} B \in M_m(\C) = Y .
\end{equation*}
Now it is elementary to verify that any fiber $ \psi^{-1}( y ) $ contains a full rank matrix $ w $ under the condition $ m \leq n $.
\end{proof}

It is easy to see that $ Y = U \otimes U \simeq M_m(\C) $ inherits the $ SL(m,\C) $-action in Eq.~\eqref{eq:action.on.matrices} 
considered by {\Djokovic}-Sekiguchi-Zhao \cite{DSZ.2003}, 
while $ X $ is isomorphic to $ \AffineGrass_m( V \oplus V^{\ast} ) $.  
Thus, an $ SL(m, \C) $-orbit $ \Osmall \subset M_n(\C) $ is lifted to 
a $ GL(n, \C) $-orbit $ \Obig = \theta( \Osmall ) \subset \AffineGrass_m( V \oplus V^{\ast} ) $.  
Moreover, the lifting map 
\begin{equation}
\theta : M_m(\C) / SL(m, \C) \to \AffineGrass_m( V \oplus V^{\ast} ) / GL(n, \C) 
\end{equation}
is injective, and preserves the closure relation and nilpotency.

As an example, we examine the simplest case, i.e., the lifting of the trivial orbit.

\begin{example}
Let $ \Obigtrivial = \theta ( \{ 0 \} ) $ be the lift of the trivial nilpotent orbit $ \Osmall = \{ 0 \} \subset Y $.  
Then we have $ \closure{\Obigtrivial} = \nullcone \GITquotient \Gsmall $.  
Let us show that $ \Obigtrivial $ is a $ GL(n, \C) $-spherical variety 
with the normal closure in $ X = \AffineGrass_m( V \oplus V^{\ast} ) $.  

To prove it, let us prepare some notations.  
We denote the set of partitions of length at most $ m $ by 
\begin{equation*}
\partition_m = \{ \alpha = ( \alpha_1, \ldots, \alpha_m ) \in \Z^m \mid 
\alpha_1 \geq \cdots \geq \alpha_m \geq 0 \} , 
\end{equation*}
which we identify with the subset of dominant weights of $ GL(m, \C) $ as usual.  
The irreducible finite dimensional representation of $ GL(m, \C) $ with highest weight $ \alpha $ is denoted by $ \tau^{(m)}_{\alpha} $, and 
its contragredient $ \tau^{(m)}_{\alpha}{}^{\ast} $ has the highest weight 
$ \alpha^{\ast} = ( - \alpha_m, - \alpha_{m - 1} , \ldots, - \alpha_1 ) $, 
hence $ \tau^{(m)}_{\alpha}{}^{\ast} = \tau^{(m)}_{\alpha^{\ast}} $ holds.  
For $ \alpha, \beta \in \partition_m $, we put 
\begin{equation*}
\alpha \composit \beta = ( \alpha , 0 , \ldots, 0, \beta^{\ast} ) \in \Z^n ,
\end{equation*}
which is a dominant weight of $ GL(n, \C) $.  

Now let us return to the proof.  
We note that as a $ GL(n, \C) \times ( GL(m, \C) \times GL(m, \C) ) $-module,
\begin{equation*}
\C[ \nullcone ] \simeq \harmonics \simeq \directsum_{\alpha, \beta \in \partition_m} \tau_{\alpha \composit \beta}^{(n)} \boxtimes ( \tau_{\alpha^{\ast}}^{(m)} \otimes \tau_{\beta^{\ast}}^{(m)} ) .
\end{equation*}
For this, see Theorem~2.5.4 in \cite{Howe.1995} for example 
(note that the action of $ GL(m, \C) $ in the second factor is the dual to that in \cite{Howe.1995}).  
Then the action of $ SL(m, \C) $ is obtained by the restriction of the action of $ GL(m, \C) \times GL(m, \C) $ to the diagonal subgroup $ \Delta SL(m, \C) $.  
Note that 
\begin{equation*}
\begin{split}
( \tau_{\alpha^{\ast}}^{(m)} \otimes \tau_{\beta^{\ast}}^{(m)} )^{SL(m, \C)} 
&= \Hom_{SL(m, \C)}( \tau_{\alpha}^{(m)} , \tau_{\beta^{\ast}}^{(m)} ) \\
&= \begin{cases}
\C & \text{if $ \alpha - \beta^{\ast} = k \trivial_m $ for some $ k \in \Z $} , \\
0  & \text{otherwise} , \\
\end{cases}
\end{split}
\end{equation*}
where $ \trivial_m = ( 1 , 1 , \ldots , 1 ) \in \Z^m $.  
Thus we have $ ( \tau_{\alpha^{\ast}}^{(m)} \otimes \tau_{\beta^{\ast}}^{(m)} )^{SL(m, \C)} \simeq \C $ if and only if 
$ \alpha_i + \beta_{m - i + 1} = k \; ( 1 \leq i \leq m ) $ for a fixed non-negative integer $ k $.
From this, we obtain the decomposition of the regular function ring $ \C[ \closure{\Obigtrivial} ] $ of the closure of the lifted orbit $ \Obigtrivial $ as a representation of $ GL(n, \C) $.  
\begin{equation}
\C[ \closure{\Obigtrivial} ] \simeq \C[ \nullcone ]^{SL(m, \C)} 
     \simeq \directsum_{ k \geq 0 } \directsum_{\alpha - \beta^{\ast} = k \trivial_m} \tau_{\alpha \composit \beta}^{(n)} \quad 
(\text{$ \alpha, \beta $ move over $ \partition_m $}) 
\end{equation}
This shows that $ \closure{\Obigtrivial} $ is a spherical variety, hence so is $ \Obigtrivial $.  
As for the normality, $ \nullcone $ is known to be normal, and as a quotient of the normal variety, $ \closure{\Obigtrivial} $ is also normal.
\end{example}

\subsection{Resolution via the action of the orthogonal and symplectic groups.}

We denote the natural representation of $ O(p, \C) $ by $ V^+ $ and that of $ Sp(2 q, \C) $ by $ V^- $.  
We put $ \Gbig = G_+ \times G_- = O(p, \C) \times Sp(2 q, \C) $ and 
$ W = ( V^+ \oplus V^- ) \otimes U $ for the natural representation $ U = \C^m $ of $ \Gsmall = SL(m, \C) $.  
We assume that $ 2 m < p, 2 q + 1$. 
By Theorem \ref{thm:tensor.product.type}, the double fibration by the affine quotient maps 
\begin{equation*}
X = W \GITquotient \Gsmall \xleftarrow{ \quad \varphi \quad } W \xrightarrow{ \quad \psi \quad } W \GITquotient \Gbig = Y
\end{equation*}
satisfies Assumptions \ref{basic.assumption.Y} and \ref{basic.assumption.X}.  
Moreover, we have 
\begin{equation*}
Y = \Sym(U) \oplus \Alt(U) \simeq M_m(\C) 
\end{equation*}
with the $ SL(m,\C) $-action in Eq.~\eqref{eq:action.on.matrices}, 
while $ X $ is isomorphic to $ \AffineGrass_m( V^+ \oplus V^- ) $.  
Thus, an $ SL(m, \C) $-orbit $ \Osmall \subset M_n(\C) $ is lifted to 
an $ O(p, \C) \times Sp(2 q, \C) $-orbit $ \Obig \subset \AffineGrass_m( V^+ \oplus V^- ) $.  
Moreover, the lifting map 
\begin{equation}
\theta : M_m(\C) / SL(m, \C) \to \AffineGrass_m( V^+ \oplus V^- ) / O(V^+) \times Sp(V-)
\end{equation}
is injective and it preserves the closure relation and maps nilpotent orbits to nilpotent orbits.  

Thus, if we put $ p = 2 q $, then the orbit space $ M_m(\C) / SL(m, \C) $ has two different embeddings to the same affine Grassmanian cone 
$ \AffineGrass_m( V^+ \oplus V^- ) $.  One is treated in this subsection, 
and the other is explained in \S \ref{subsec:resolution.by.contraction}.

%\bibliographystyle{plain}
%\bibliography{tuebingen2003proc}

\end{document}